\documentclass[10pt,confmode]{IREP2013}

\usepackage{multicol}
\usepackage[pdftex]{graphicx}
\usepackage{amssymb,amsfonts,amsmath}
\usepackage{eucal,bm}
\usepackage{gensymb,subfigure}
\usepackage{color}
\usepackage{booktabs}
\usepackage{url}

\newcommand{\cE}{{\cal E}}

\def\cB{{\cal B}}

\newcommand{\R}{\mathbb{R}}

\newcommand{\beq}{\begin{equation}}
\newcommand{\eeq}{\end{equation}}
\newcommand{\beqnr}{\begin{eqnarray}}
\newcommand{\eeqnr}{\end{eqnarray}}

\newcommand{\benum}{\begin{enumerate}}
\newcommand{\eenum}{\end{enumerate}}

\newcommand{\cL}{{\cal L}}

\newtheorem{DE}{Definition}[section]

\newtheorem{LE}[DE]{Lemma}

\newcommand{\edges}{{\cal E}}
\newcommand{\nodes}{{\cal V}}
\newcommand{\generation}{p}
\newcommand{\demand}{d}
\newcommand{\wind}{w}
\newcommand{\phase}{\theta}
\newcommand{\auxphase}{\vartheta}
\newcommand{\voltage}{v}
\newcommand{\maxgeneration}{p^{max}}
\newcommand{\mingeneration}{p^{min}}
\newcommand{\generators}{\mbox{\cal G}}
\newcommand{\demands}{\mbox{\cal D}}
\newcommand{\windgenerators}{\mbox{\cal W}}
\newcommand{\susceptance}{\beta}
\newcommand{\capacity}{\bar{p}}
\newcommand{\cost}{f}
\newcommand{\lossflow}{\rho}
\newcommand{\arcsinflow}{\psi}
\newcommand{\lagrange}{\theta}
\newcommand{\droop}{\alpha}
\newcommand{\generationsetpoint}{p}

\newif\ifcomments \commentstrue
\newcommand\checkspace{\ifvmode\vspace{80pt}\penalty0\vspace{-80pt}\leavevmode\fi}
\newcommand\defedit[3]{\newcommand#1[1]{\ifcomments{\color{#2}{\checkspace
                                                    \ifinner{{\color{#2}{[#3]: }}}\else\marginpar{{\color{#2}{#3}}}\fi
                                                   \bf{##1}}}\fi}}
\defedit{\rbc}{red}{RB}
\defedit{\dbc}{blue}{DB}
\defedit{\mcc}{green}{MC}

\begin{document}
\onecolumn
\title{\Large\textbf{Synchronization-Aware and Algorithm-Efficient\\ Chance Constrained Optimal Power Flow}}

\vspace{.2cm}
\author{\begin{tabular}{ccc}
{\bf Russell Bent} & {\bf Daniel Bienstock} & {\bf Michael Chertkov}\\
\begin{tabular}{c} D-4 of Defense Systems\\ and Analysis Division,\\
Los Alamos NL,\\ NM 87545 USA\\ ({\tt rbent@lanl.gov})
\end{tabular}
 &
\begin{tabular}{c}
\\ Departments of Industrial\\ Engineering \& Operations Research,\\ Applied Physics and Mathematics,\\
Columbia University, 500 West 120th St.\\
New York, NY 10027 USA\\ ({\tt dano@columbia.edu})
\end{tabular}
&
\begin{tabular}{c}
\\ T-4 of Theoretical Division\\ and Center for Nonlinear Studies,\\ Los Alamos NL,\\
and New Mexico Consortium,\\ NM 87545(4) USA\\ ({\tt chertkov@lanl.gov})
\end{tabular}
 \end{tabular}}

\maketitle

\begin{multicols}{2}

\section{Abstract}

One of the most common control decisions faced by power system operators is the question of how to dispatch generation to meet demand for power.  This is a complex optimization problem that includes many nonlinear, non convex constraints as well as inherent uncertainties about future demand for power and available generation.  In this paper we develop convex formulations to appropriately model crucial classes of nonlinearities and stochastic effects.  We focus on solving a nonlinear optimal power flow (OPF) problem that includes loss of synchrony constraints and  models wind-farm caused fluctuations. In particular, we develop  (a) a convex formulation of the deterministic phase-difference nonlinear Optimum Power Flow (OPF) problem; and (b) a probabilistic chance constrained OPF for angular stability, thermal overloads and generation limits that is computationally tractable.


\section{I. Introduction}

Generation re-dispatch is a routine operation for adjusting the output of flexible generators to meet the needs for electric power in transmission systems.  Redispatch is performed periodically based on predictions about the state of system over a planning horizon, i.e. demand for power over the next time period (usually 15 minutes to an hour). These operations rely on computed solutions to variations of the Optimal Power Flow (OPF) problem.  Typical implementations of OPF minimize the (convex) cost of generation, subject to linearized power flow balance constraints, thermal capacity constraints (on power lines), ramping constraints, security constraints, etc.  In  recent work \cite{12BCH}, we have developed a chance constraint OPF model (CC-OPF) that generalizes the standard OPF to include uncertainties from fluctuations in the output of wind farms.  This CC-OPF manages risk in a principled fashion by allowing physical constraints to be violated (such as line capacity constraints) with small and controlled probability. Furthermore, our work has shown how to construct a robust version of the CC-OPF -- optimal within the set of power flows valid for a range of parameters characterizing the probability distribution of the wind output. In practice, unconstrained solutions to an OPF produce a relatively small number of overloads that violate the chance constraints, making this problem well-suited for the cutting plane algorithm we developed in \cite{12BCH}.  This algorithm can solve large instances of the CC-OPF, such as the 2746 bus Polish network in 20 seconds using a desktop computer.

As a standard simplification practice, OPF models use the linearized DC equations for modeling power flow (PF) physics.  The DC approximation contains a number of key assumptions (see e.g. \cite{gmo07power}): (a) voltage is constant (fixed) at all nodes of the network; (b) thermal (resistive) losses are negligibly small (ignored); and (c) the phase difference over any line of the network is small, $|\theta_i-\theta_j|\ll 1$, and thus the power flows are linearized in the phases.

In this paper, we extend the standard (deterministic) OPF approach and the CC-OPF approach of \cite{12BCH} to account for some of the nonlinear effects in PF physics. Specifically, in describing the PFs we do not include assumption (c).  In other words, we still assume that voltage is constant and lines are lossless; however, we include the correct nonlinear phase difference model in the OPF computation.

First, for the AC (phase-nonlinear) deterministic PF we develop two complementary approaches for solving the deterministic AC-OPF. In the first approach, we utilize the results of \cite{12DCB} and show how the  nonlinear synchronization constraint can effectively be modeled using a linear set of inequalities.  The computational complexity of the resulting model is equivalent to the standard convex DC-OPF algorithm. In the second approach, we adapt and extend the convex formulation of the AC PF in \cite{12BV} to derive an exact, nonlinear but still convex (thus computable efficiently) AC-OPF.

Second, we generalize the CC-OPF model developed by \cite{12BCH} to incorporate the synchronization condition developed in \cite{12DCB}. We show that
the resulting ``sync-aware'' CC-OPF is reducible to a convex optimization problem with the same complexity as the CC-OPF model of \cite{12BCH}. We finally develop a computationally very efficient algorithm for this problem that is capable of solving optimal generation re-dispatch for networks with thousands of nodes  in seconds on a standard computer.

Our model includes chance constraints associated with thermal limits and
with the synchronization condition.  The latter express probabilistically
the requirement that fluctuations in wind resources may cause the system to lose synchrony for only very short periods of time (second or less).
In most dense systems (like the Eastern interconnect in US or European grids) the synchronization condition is typically less constraining than the respective conditions guarding against thermal overloads of lines.  In sparse systems with long lines this situation is often reversed; for example transmission grids of Russia, Australia or some of Midwest and Southwest states in US are loss-of-synchrony prone. In our experiments we include  constraints of both types to ensure both limits are addressed. We illustrate this competition on examples of small, moderate and large sample grids.

Within the power engineering literature there have been a number of recent papers developing chance constrained versions of various problems. First \cite{zl11chance} developed a chance constrained version of DC-OPF that focuses on variations in demand.  Generation is dispatched to ensure that the system's physical constraints are violated with low probability (the chance constraints).  A local search algorithm is developed to solve the problem.  \cite{12BCH} develops a convex chance constrained OPF for problems with uncertain renewable energy fluctuations using a closed form expression of the fluctuations. Chance constraints have also been used in longer time scale problems such as unit commitment \cite{omn04unit, wsl08unit,zshy11unit} as well as expansion planning \cite{yw05chance, ycwc09tnep} to account for uncertainty in renewable energy generation.  Finally, \cite{hgrz12voltage} develops a chance constraint model for voltage control.

Optimal power flow problems have been the subject of a large amount of research \cite{744492,744495}. Most closely related to this paper is recent work that has developed approaches for solving versions of the AC OPF using conic and convex relaxations.  As examples, there interior point methods such as in \cite{jabr08} and \cite{ll11opf} that develops a dual approximation to the ACOPF which is convex semi-definite, along with many others. There has also been literature on adding constraints to the DC model in order to generate solutions that are close to their AC counterpart \cite{cbvh12smart,cvhb12approx}.

The remainder of this paper is organized is follows.  Section II describes notation we use for describing power networks.  Section III describes how turn the AC OPF into a convex problem.  Section IV adds the chance constraints.  Section V describes our algorithm and Section VI some experimental results.  We conclude and discuss path forward with Section VII.

\section{II. Model}

A power network is defined by a set of nodes (buses) and edges (power lines), referred to here respectively as $\nodes$ and $\edges$.  For each node $\i \in \nodes$, we use $\generation_i, \demand_i, \wind_i, \theta_i$, and $\voltage_i$ to denote the conventional generation, power consumption, wind generation, phase angle, and voltage at $i$.  Similarly, $\mingeneration_i$ and $\maxgeneration_i$ is used to denote the minimum and maximum amount of conventional generation at $i$. We use the sets (possibly overlapping)
$\generators = i \in \nodes : \maxgeneration_i \neq 0$,  $\demands = i \in \nodes : \demand_i \neq 0$,  $\windgenerators = i \in \nodes : \wind_i \neq 0$
to denote sets of nodes that, respectively,  generate power at traditional power plants (controllably), consume power and generate power from fluctuating/uncertain (wind) sources. Finally the cost of generating power for $i \in \generators$ is denoted by $\cost_i$ and $\droop_i$ is used to model affine response of controllable generation (through primary and secondary control of frequency at the generators.

An (undirected) edge $ij \in \edges$ is defined by the nodes $i$ and $j$ it connects. The terms $\susceptance_{ij}$ and $\capacity_{ij}$ are then used to define the susceptance and thermal capacity of $ij$. In many places, we assume, as
customary that an arbitrary orientation of the edges has been chosen, yielding a directed
graph; any arc of graph is denoted by $(i, j)$ if it has node $i$ as the ``tail'' or ``from'' node and node $j$ as the ``head'' or ``to'' node.  For convenience, given an arc $k$ we will write $\beta_k$ and $\capacity_k$ for the susceptance (resp., capacity) of the (undirected) edge corresponding to arc $k$.
Let $\cL$ denote the set of arcs, and construct the node-arc incidence matrix of this directed graph as follows:
for any node $i$ and arc $k = (u,v)$ we have $a_{ik} = 1$ if $i = u$,
$a_{ik} = -1$ if $i = v$ and $0$ otherwise.  Arc and edge notation is used interchangeably throughout the text, in order to keep the formulations compact.

\section{III. AC OPF Convexification}

Our paper uses the real power, thermally-constrained, loseless and voltage-constrained  optimal power flow (OPF)  problem as our starting point:

\begin{alignat}{2}
\min\limits_{\generation,\phase} \ &\sum\limits_{i\in \generators} \cost_i (\generation_i) = \min\limits_{\generation,\phase}\cost(\generation) \label{AC_OPF} \\
\text{s.t.} &  \  \!\sum\limits_{j: ij\in \edges} \!\!\sin(\phase_i-\phase_j)\susceptance_{ij}\!=\! \generation_i\! -\! \demand_i\! +\!\wind_i   & \ \ \ \forall_{ i\in \nodes} \label{PF} \\
 &  |\sin(\phase_i-\phase_j)\susceptance_{ij}| \le \capacity_{ij} & \ \ \forall_{ij\in \edges} \label{T1} \\
 &  \mingeneration_i \leq \generation_i\leq \maxgeneration_i  & \ \  \forall_{i\in \nodes} \label{G1}
  \end{alignat}

Eq.~(\ref{AC_OPF}) minimizes the cost of dispatching conventional generation (usually convex, quadratic) to meet $\demand$ and $\wind$. Eqs.~(\ref{PF}) states the  Power Flow (PF) constraints under the assumption of uniformly maintained voltage\footnote{In this ACOPF, voltages are modeled as per unit and assumed to be 1, though there is nothing in the model that prevents us from using fixed voltages other than 1.} and of purely inductive (zero resistivity) lines\footnote{The generation and consumption of the ACOPF is zero sum, i.e. $\sum_{i\in\nodes}(\generation_i-\demand_i+\wind_i)=0$,  reflecting the lossless nature of the network.}. Eqs.~(\ref{T1}) and Eqs.~(\ref{G1}) introduce constraints on the line flows and conventional generation. Finally, we note that this OPF model is exactly the same as the traditional DCOPF, except for the more accurate sine in (\ref{PF}). We now discuss two ways of transforming this problem into a convex problem.

\subsection{Convex Optimization of Voltage Uniform and Lossless Power Flows}

We first describe a provably correct approach for directly transforming the above OPF problem into a convex optimization problem. This transformation is based on a reformulation of Eqs.~(\ref{PF}) first discussed in \cite{12BV}.

Our formulation is:
\begin{eqnarray}
&& \min\limits_{\lossflow} \  \sum_{k \in \cL} \susceptance_{k} \arcsinflow(\lossflow_{k})  \label{IndLoss} \\
&& \text{s.t.} \ \  \sum_{k \in \cL} \susceptance_{k} a_{k} \lossflow_{k}=\generation_i-\demand_i+\wind_i   \ \ \ \forall i\in \nodes  \label{conserv} \\
 &&  |\lossflow_{k}| < \min\{1,\capacity_{k}\}.  \ \ \forall_k \in \cL \label{rho_less_1}
  \end{eqnarray}
\noindent Here for $|x| < 1$,
$$\arcsinflow(x)\doteq\int_{-1}^x \arcsin(y)\, dy,$$
a convex function of $x$ since $\arcsin(x)$ is increasing for $x \in [-1, 1]$. Interestingly, in this formulation, if the optimal solution of Eq.~(\ref{IndLoss}) occurs on the boundary of Eq.~(\ref{rho_less_1}), then there exists no feasible solution to Eqs.~(\ref{PF}). In other words,  the grid cannot be synchronized.

When optimization problem (\ref{IndLoss}-\ref{rho_less_1}) is well-defined, that is to say, it has an optimal solution which satisfies the strict Eqs.~(\ref{rho_less_1}), we can verify that Eq.~(\ref{IndLoss}) yields Eqs.~(\ref{PF}).  To see this, we use convexity of the objective (\ref{IndLoss}); denoting by $\theta_i$ the optimal Lagrange multiplier for constraint (\ref{conserv}) we obtain that for every arc $k = (i, j) \in \cL$,
\begin{eqnarray}
\quad \arcsin(\lossflow_{k})=\lagrange_i-\lagrange_j. \ \label{rho-theta}
\end{eqnarray}
The combination of this result with Eqs.~(\ref{conserv}) renders Eqs.~(\ref{PF}) \footnote{See also Appendix A for a brief discussion of the relation between the optimization formulation just discussed and the so-called energy function approach.}. Unlike many convex transformations, there is a physical meaning to the optimization in Eq.~(\ref{IndLoss}). It is twice the reactive power losses in all the lines of the network. Finally, the constraints (\ref{conserv}) expresses flow conservation at any node of the network.

\subsection{Convex Optimization of Voltage-Uniform and Lossless-Optimal Power Flows}

We now combine the derivation in the prior section with optimal power flow.  Let $C$ be lower bound on the optimal OPF cost\footnote{This can be obtained, for example, by solving a DCOPF problem.}. Let $0 < \epsilon < 1$ be a tolerance, $\beta_{max} = \max_{k} \beta_{k}$, $D = \frac{C \epsilon}{\pi \beta_{max}}$ and $\phi = (-m \log \epsilon)^{-1} \beta_{max}$, where $m$ = number of lines. For an arc $k$, let its {\em effective capacity} be $u_k = \min\{1,\capacity_{k}/\beta_{k}\}$. Finally, let the {\em barrier} function $B \, : \, \R_+ \rightarrow \R_+$ be defined by
$$ B(t) \ = \ - \log t.$$
Consider the optimization problem given by
\begin{eqnarray}
&& \min \ f(p) \ + \ D \sum_{ k \in \cL} \susceptance_k
 \left[ \arcsinflow(\lossflow_{k})  \ + \ \phi \, B(\delta_{k}) \right] \label{newobj} \\
&& \text{s.t.} \ \sum_{k \in\cL} \susceptance_{k}a_{ik} \lossflow_{k}=\generation_i-\demand_i+\wind_i   \ \ \ \forall i\in \nodes  \label{conserv2} \\
&& |\lossflow_{k}|\ + \  u_k \delta_{k} \ \le \ u_k  \label{rho_less_1_2}\\
&& \delta_{k} \ \ge \ 0, \ \ \forall k \in \cL. \label{newlast}
\end{eqnarray}
We have:
\begin{LE}  (a) Suppose $(p^*, \theta^*)$ is the optimizer for the sync-constrained OPF problem. Suppose, further, that for each arc $k = (i,j)$ we have
\begin{eqnarray}
&& | \sin(\theta^*_i - \theta^*_j) | \ \le \ (1 - \epsilon) u_k.
\label{slacksine} \end{eqnarray}
Then, defining $\rho^*_{ij} = \sin(\theta^*_i - \theta^*_j)$ for each arc
$(i,j)$, we obtain that $p^*, \rho^*$ is a feasible solution to problem
(\ref{newobj})-(\ref{newlast}) with cost at most
$$ (1 + 2 \epsilon) f(p^*).$$
(b) Conversely, let $(\hat p, \hat \rho, \hat \delta)$ be an optimal solution to problem (\ref{newobj})-(\ref{newlast}) with objective value $\hat K$, say. Then, with $\hat \theta$ obtained by solving problem (\ref{IndLoss})-(\ref{rho_less_1}) we have that $(\hat p, \hat \theta)$ is feasible for the sync-constrained OPF problem and has cost $f(\hat p) \le \hat K$. (c) Let $(\hat p, \hat \rho, \hat \delta)$ be as in (b), and suppose that additionally for every arc $k$ we have
\begin{eqnarray}
&& |\hat \rho_k| \le u_k (1 - \tilde \epsilon), \label{separation}
\end{eqnarray}
for some value $0 < \tilde \epsilon < 1$.  Let $D \hat \theta$ be the optimal dual variables for constraints (\ref{conserv2}).  Then for every line $k = (i,j)$ we have:
\begin{eqnarray}
&& \hat \rho_k \ = \ \sin(\hat \theta_i - \hat \theta_j \ + \ \eta_k) \label{approxsine}
\end{eqnarray}
where $ | \eta_k | \le (m \, u_k \, \tilde \epsilon \log (1/\epsilon) )^{-1}$.
\end{LE}
\noindent {\em Proof.}  (a) Define, for each arc $k = (i,j)$, $\rho^*_{k} = \sin(\theta^*_i - \theta^*_j)$.  Since (\ref{slacksine}) holds the values
$$\delta^*_{k} = 1 - \frac{1}{u_k} \, |\lossflow^*_{k}|$$
satisfy
$$ \epsilon \le \delta^*_{k} \le 1,$$
and the vector $p^*$, $\rho^*$, $\delta^*$ is feasible for (\ref{newobj})-(\ref{newlast}).  It follows that for any arc $k$,  $-B(\delta^*_{k}) \le \log(1/\epsilon)$, and so the total contribution from all the terms $B(\delta^*_{k})$ in (\ref{newobj}) is at most $D \le \frac{\epsilon}{\pi} C < \epsilon f(p^*)$.  Likewise, the total contribution from the terms involving the $\arcsinflow(\lossflow^*_{k})$ is also at most $\epsilon f(p^*)$.

\noindent (b) Follows directly from the structure of the objective function (\ref{newobj}).

\noindent (c) The first-order optimality conditions for each line $k = (i,j)$
yield
$$ \hat \theta_i - \theta_j = \arcsinflow(\hat \lossflow_{k}) \ - \ \frac{\phi}{u_k} \frac{1}{1 - \hat \lossflow_k/u_k}.$$
Since $ 1 - \hat \lossflow_k/u_k \ge \tilde \epsilon$ the result follows.
\QED

{\bf Remarks:} Condition (\ref{separation}) amounts to a {\em separation} condition: we find a solution where each flow quantity $| \hat \rho_k|$ is sufficiently separated from the effective capacity $u_k$.  For small but fixed $\epsilon$, larger $\tilde \epsilon$, and $m$ large, the quantity $\eta_k$ in (\ref{approxsine}) is small in absolute value.  Thus condition (\ref{approxsine}) indicates that the flows and (scaled) dual variables computed in the optimization problem (\ref{newobj}-\ref{newlast}) yield an approximately feasible solution to the OPF problem (\ref{AC_OPF}-\ref{G1}).  Also note that the barrier terms in the objective function (\ref{newobj}) tend to $+\infty$ as $|\rho_k| \rightarrow 1$ for some $k$; from a practical perspective this (together with the above Lemmas) imply that the value of problem (\ref{newobj}-\ref{newlast}) becomes ``large'' precisely when the ACOPF problem (\ref{AC_OPF}-\ref{G1}) becomes infeasible.  Finally, the overall approach is easily adapted to the case where we want to impose general constraints of the form $| \theta_i - \theta_j | \le \gamma_{ij}$ where the $\gamma_{ij} \le \pi/2$ are input data.

\subsection{Synchronous Constrained OPF}

The approach in the preceding section provides a provably correct modification of the ACOPF problem (\ref{AC_OPF})-(\ref{G1}) into a convex optimization problem.  Even though the problem is indeed convex, it may not necessarily be trivially solvable \footnote{In future work we plan to test the computational feasibility of solving problem (\ref{newobj})-(\ref{newlast}).}. Moreover we do not have a way to extend the approach so as to handle additional constraints on the power flow.  In this section we consider a different (admittedly heuristic) approach that also renders a convex optimization problem which in particular can handle side-constraints; specifically stochastic constraints.

Returning to the original OPF problem (\ref{AC_OPF1})-(\ref{G1}), the main challenge with this model is that  Eqs.~(\ref{PF},\ref{T1}) are nonlinear. One approach for addressing this difficulty replaces Eqs.~(\ref{PF},\ref{T1}) with equivalent linear, expressions. To achieve this result we modify Eq.~(\ref{AC_OPF}) according to the results of \cite{12DCB}.  Reference \cite{12DCB} discovered the following condition for lossless and voltage-maintained power flow models:\\

\underline{\bf Synchronous (Sync) condition \cite{12DCB}}

For a fixed $p,d,w$, consider the following system of equations:
\begin{eqnarray}
&& \ |\auxphase_i-\auxphase_j|< \min\left(\capacity_{ij}/\susceptance_{ij},1\right) \ \ \ \forall ij\in \edges
\label{T2}\\
&& \ \sum_{j\sim i}(\auxphase_i-\auxphase_j)\susceptance_{ij}=p_i-d_i+w_i \ \ \ \ \forall i\in \nodes.
\label{vartheta}
\end{eqnarray}
Given a solution to this system, suppose we can choose for each node $i$ a value $\phase_i$ such that for each edge $ij$, $\sin(\phase_i - \phase_j) = \auxphase_i - \auxphase_j$.  Then the $\phase_i$ solve Eqs.~(\ref{PF})- (\ref{G1}).

The work in \cite{12DCB} suggests,  via analytical and statistical methods, that in many cases this change leads to an accurate computation of line flows, though possibly not of actual phase angles. For example, the change of variables is exact on trees. We thus obtain a formulation that is reminiscent of the traditional DC formulation, however accounting for the nonlinearities in Eqs. (\ref{PF},\ref{T1}).

Thus we are lead to replace  Eqs.~(\ref{PF},\ref{T1}) with  Eqs.~(\ref{vartheta},\ref{T2}) in the AC OPF (\ref{AC_OPF}), or in other words, to substitute Eq.~(\ref{AC_OPF}) by the following convex optimization problem
\begin{equation}
\min_{p,\vartheta} f(p),\mbox{    s.t.   }
\mbox{Eqs.~(\ref{G1},\ref{T2},\ref{vartheta})}. \label{AC_OPF1}
\end{equation}
Remarkably, Eq.~(\ref{AC_OPF1}) is identical to the DCOPF with an extra constraint on the flows across an edge. We therefore refer to this OPF as  the sync-constrained OPF (SCOPF).

\section{IV. Chance Constrained OPF}

In the previous sections we assumed that the load ${\bf \demand}$ and uncontrollable generation ${\bf \wind}$ are known apriori and do not significantly or unpredictably change over a dispatch planning horizon;  this assumption is justified in practice. In contrast, wind generation uncertainty  is  significant and cannot be ignored \cite{08BDL,CIGRE09,11HM}.  Wind fluctuations are intrinsically stochastic phenomena; even when the wind forecast  is known,  the forecast only expresses information on the underlying probability distribution but not on precise wind generation.

To overcome this limitation \cite{12BCH} developed a stochastic model of wind and other uncertain resources for OPF problems. More formally,  the  Chance Constrained (CC) OPF with fluctuating wind resources is stated as

\begin{alignat}{2}
\min\limits_{\generationsetpoint,\droop} \ & \mathbb{E}_w\left[ f(\generation)\right] \label{CC-schem}  \\
\text{s.t.} &  \     \mbox{Prob (CON violation)} < \varepsilon_{\mbox{\tiny CON}}       & \ \ \ \forall \ {\mbox{CON}} \label{CC}
\end{alignat}

where $\mathbb{E}_w\left[\cdots\right]$ is used to denote the expectation over a probability density function (PDF) for the wind, $\rho(\wind)$, s.t. $\int \demand \wind \rho(\wind)=1$. In (\ref{CC}) CON represents an OPF constraint such as a (thermal) line limit or a generator output limit; $\varepsilon_{\mbox{\tiny CON}}$ is a small number (additional parameter) controlling the probability that the constraint is violated. In (\ref{CC}) $\generationsetpoint,\droop$ are vectors describing the generation set points and affine (droop) coefficients which will be discussed, in detail, later in the text.

The CC-OPF computation that serves as our inspiration was discussed in \cite{12BCH} in the context of the linear DC power flow model.  This approach developed chance constraints for thermal constraints on power lines and chance constraints for bounds on generation. Within the DC model all input configurations (of $\demand-\wind$) are guaranteed to have a solution satisfying power transmission constraints (but possibly not thermal line or generator limit constraints).  However,  in nonlinear models such as those discussed in the previous section, the power flow Eqs.~(\ref{PF}) may have no solution -- equivalently the system may be out of sync. More formally, this CC-OPF version of Eq.~(\ref{CC}) is stated as

\begin{alignat}{2}
\min\limits_{p,\alpha} \ &  \mathbb{E}_{\bm{\wind}}\left[f(\generationsetpoint - (e^T \bm{\wind})\droop ) \,  \right]  \label{CC-OPF} \\
\text{s.t.} &  \   \sum\limits_{i\in\generators} \droop_i  =  1,\quad  \droop \ge 0,  \quad \generationsetpoint \ge 0    \label{conic-first} \\
& \quad \mbox{Prob} \left( \mbox{PF Eqs. are not feasible}\right) <  \varepsilon &\  \  \label{sync}\\
& \quad \mbox{Prob} \left( \susceptance_{ij}|\sin(\phase_i-\phase_j)| > \capacity_{ij}\right) <  \epsilon_{ij} &\  \ \forall_{i,j \in \edges} \label{longchance} \\
& \quad \mbox{Prob} \left(\generationsetpoint_g - (e^T\bm{\wind}) \droop_i > \maxgeneration_g \right) < \epsilon_g &\  \forall_{g \in \generators} \label{chance_gen_min} \\
& \quad \mbox{Prob} \left(\generationsetpoint_g - (e^T\bm{\wind}) \droop_i < \mingeneration_g \right)  < \epsilon_g &\  \forall_{g \in \generators} \label{chance_gen_max}
\end{alignat}
where $\mu=\mathbb{E}_\wind[\wind]$ is the mean wind, ${\bm \wind}=\wind-\mu$ is the zero mean fluctuating component of the wind, and $\generationsetpoint-(e^T \bm{\wind})\droop$ is the vector of controllable generation according to the affine model of the automatic (primary and secondary) generation control. PF entering Eqs.~(\ref{longchance},\ref{chance_gen_min},\ref{chance_gen_max}) assume the following affine response version of Eqs.~(\ref{PF})
\begin{alignat}{2}
& \  \sum\limits_{j: (i,j)\in \edges} \sin(\phase_i-\phase_j)\susceptance_{ij} \nonumber \\
&=\generationsetpoint_i -(e^T\bm{\wind})\droop_i- d_i +\mu_i+\wind_i & \ \ \ \forall_{i\in \nodes}.
\label{PF1}
\end{alignat}

In  Eq.~(\ref{longchance}), there are a number of ways that $\epsilon_{ij}$ can be interpreted.  For the purposes of this paper, we interpret $\epsilon_{ij}$ as the fraction of time that a line's flow exceeds a critical value, such as its thermal limits.  For example, setting $\epsilon_{ij} = \frac{1}{60}$ for thermal limits is equivalent to stating that a line may be overloaded for at most 1 minute during a 1 hour period (a typical planning horizon for dispatch decisions).

The $\varepsilon$ in the probabilistic sync condition (\ref{sync})  is interpreted  similarly to $\epsilon_{ij}$ entering the probabilistic thermal constraints
(\ref{longchance}). It is the fraction of time a system can tolerate loss of synchrony. However, the  quantitative difference between the sync and thermal constraints (in actual values) is significant.  In practice, loss of synchrony cannot last for more than a second, or even a fraction of a second, before electro-mechanical instability develops in the system. Therefore, the $\varepsilon$ associated with loss of synchrony is typically significantly ($\sim$two orders of magnitude) smaller than, $\varepsilon_{ij}$.  In summary,  constraint (\ref{sync}) may be active (optimization limiting), for systems with long lines, in spite of the fact that the loss of synchrony threshold is often more stringent than the thermal constraint, i.e.  when $\bar{p}_{ij}/\susceptance_{ij}<1$, which holds for typical lines, $\varepsilon$ smaller than $\epsilon_{ij}$ makes Eq.~(\ref{sync}) potentially more restrictive than Eq.~(\ref{longchance}).

The optimization problem described in (\ref{CC-OPF}) is non-convex, and generally difficult to solve. However, by using the modifications, conjectures and simplifications developed for the SCOPF, the problem is transformed into a formulation that is computationally tractable.

\subsection{Sync CC-OPF}

Based on the synchronous constraint developed in the previous section, we now develop a chance constrained SCOPF problem.  As in deriving SCOPF (\ref{AC_OPF1}), we replace PF Eqs.~(\ref{PF1}) with

\begin{eqnarray}
& & \sum\limits_{j: (i,j)\in{\cal E}} (\vartheta_i-\vartheta_j)\beta_{ij}\nonumber\\
& & =\bar{p}_i -(e^T\bm{\wind})\alpha_i- d_i +\mu_i+\wind_i   \ \ \ \forall_{ i\in\nodes}. \label{lin-PF}
\end{eqnarray}

Likewise, we replace Eq.~(\ref{longchance}) with (\ref{lin-PF}) and the chance constraint

\begin{equation}
 \mbox{Prob} \left(|\auxphase_i-\auxphase_j| >\bar{p}_{ij}/\susceptance_{ij}\right)  <  \epsilon_{ij} \ \ \forall_{i,j \in \edges} \label{longchance-lin}.
\end{equation}

Moreover,  following Eqs.~(\ref{T2},\ref{vartheta}) we substitute the sync feasibility chance constraint (\ref{sync}) with

\begin{equation}
\ \mbox{Prob}(|\auxphase_i-\auxphase_j|\geq 1) < \varepsilon \ \ \forall_{i,j \in \edges}. \label{sync-lin}
\end{equation}

As in SCOPF, the $\auxphase$ variables are auxiliary and are not directly related to the actual phase angles $\phase$.  The phase angles in a power flow solution resulting for a particular realization of wind are calculated using Eq.~(\ref{PF1});  this is equivalent to solving convex optimization problem (\ref{IndLoss}).
Also, once again Eq.~(\ref{CC-OPF}) only relates to the linear system of equations defined on $\auxphase$,  and in this regards it is almost equivalent to the DC CC-OPF discussed in \cite{12BCH}. The only difference is the addition of the synchrony constraint (\ref{sync-lin}).  Importantly, this new constraint  does not add complexity to the original DC CC-OPF of \cite{12BCH}.  Eq.~(\ref{sync-lin}) is a re-parameterized version of Eq.~(\ref{longchance-lin}).

In order to define the PDF we use the model of \cite{12BCH} which assumes that each component of $\wind$ is an independent, zero-mean Gaussian, i.e.

\begin{eqnarray}
\rho(\wind)=\left(\prod_{i\in{\cal W}}(2\pi\sigma_i)^{-1/2}\right)\exp\left(-\sum_{i\in{\cal W}}\frac{ \wind_i^2}{2 \sigma_i^2}\right),
\label{Gauss}
\end{eqnarray}

This allows us to evaluate all the Prob statements in Eq.~(\ref{CC-OPF}) explicitly, ie. \footnote{As in \cite{12BCH}, we do not use $\mbox{Prob}(|x|>y)<\epsilon$ directly. Instead we use
$\mbox{Prob}(x>y)<\epsilon \cup \mbox{Prob}(-x>y)<\epsilon$ which allows us to derive the convex/conic generalization described in  \cite{12BCH}}

\begin{alignat}{2}
\min\limits_{\droop,\generationsetpoint} \ & \sum_{i\in \generators} \left\{c_{i1}\left(\generationsetpoint_i^2+(\sum_{j\in \windgenerators} \sigma_j^2)
\alpha_i^2\right)+c_{i2}\generationsetpoint_i+c_{i3}\right\} \,    \label{CC-OPF-lin1} \\
\text{s.t.} &  \   \sum\limits_{i\in\generators} \droop_i  =  1,\quad  \droop \ge 0,  \quad \generationsetpoint \ge 0  \\
& |\bar{\phase}_i-\bar{\phase}_j|\leq 1-\nonumber \\
& \ \eta(\varepsilon)\left[\susceptance_{ij}^2\sum_{k\in \windgenerators} \sigma_k^2(\pi_{ik}-\pi_{jk}-\delta_i+\delta_j)^2\right]^{1/2} & \  \ \forall_{i,j \in \edges} \label{sync-lin1} \\
& \susceptance_{ij}|\bar{\phase}_i-\bar{\phase}_j|\leq \capacity_{ij}-\nonumber\\
&  \eta(\epsilon_{ij})\left[\beta_{ij}^2\sum_{k\in \windgenerators} \sigma_k^2(\pi_{ik}-\pi_{jk}-\delta_i+\delta_j)^2\right]^{1/2} & \  \ \forall_{i,j \in \edges} \label{longchance-lin1} \\
&\ \mingeneration+\eta(\epsilon_g)\left(\sum_{k\in{\cal W}}\sigma_k^2\right)^{1/2}\leq  \generationsetpoint_g &\  \forall_{g \in \generators} \label{cc-gen-min} \\
&\   \generationsetpoint_g \le\maxgeneration-\eta(\epsilon_g)\left(\sum_{k\in{\cal W}}\sigma_k^2\right)^{1/2} &\  \forall_{g \in \generators} \label{cc-gen-max}
\end{alignat}

In this formulation it is assumed that the objective function is convex-quadratic.
The term $\eta(x)$ is defined implicitly b $x=(1-\mbox{erf}(\eta(x)/\sqrt{2}))/2$.  Also
$\delta=\cB \droop$ and $\bar \phase=\breve \cB (\generationsetpoint  + \mu - d)$ \cite{12BCH}.  These terms are derived from  the $n\times n$ $\susceptance$- weighted graph Laplacian and its $(n-1) \times (n-1)$ submatrix $\hat B$ counterpart, obtained by removing row and column $n$. Finally, we also use

\begin{eqnarray}
\breve \cB & = & \left( \begin{array}{c c}
\hat B^{-1} & 0\\
0 & 0 \end{array} \right),\nonumber\\
\forall i,j\in V:&\quad & B_{ij}=\left\{\begin{array}{cc}
    -\beta_{ij}, & (i,j)\in\cE\\
    \sum_{k; (k,j)\in\cE}\beta_{kj},& i=j\\
    0,& \mbox{otherwise}\end{array}\right..
    \nonumber
\end{eqnarray}

The algorithm for solving this convex optimization problem efficiently is discussed below in Section V.

It is important to note  that even though problem (\ref{CC-OPF-lin1})-(\ref{cc-gen-max}) is convex and (as discussed below) can be efficiently solved to optimality,  its derivation relays on the extended synchrony conjecture. An alternative approach assumes that the probabilities of the chance constraints are small enough to allow the probabilities to be estimated from Large Deviation (LD) form. Some preliminary discussions of LD based approaches appear in the Appendix B.

\section{V. Algorithm}

The objective and constraints in Eq.(\ref{CC-OPF-lin1}) are convex, however since the size of the problem is large (typical practical models of transmission grids may contains thousands of nodes) it is advantageous to present the problem in a format amendable to efficient computations.  A major obstacle for efficiency is the large number of nonlinear constraints (which are nevertheless convex). To simplify computations we employed the following set of algorithmic enhancements.

First of all,  we combine the thermal Eqs.~(\ref{longchance-lin1}) and sync constraints Eqs.~(\ref{sync-lin1}), replacing these by $\forall (i,j)\in{\cal E}$:
\begin{eqnarray}
&& \left[\sum_{k\in{\cal W}}
\sigma_k^2(\pi_{ik}-\pi_{jk}-\delta_i+\delta_j)^2\right]^{1/2}\leq s_{ij},
\label{cc-nonlinear}\\
&& |\bar{\theta}_i-\bar{\theta}_j|-\eta(\epsilon_{ij}) s_{ij}\leq \capacity_{ij}/\beta_{ij},
\label{cc-lin1}\\
&& |\bar{\theta}_i-\bar{\theta}_j|-\eta(\varepsilon) s_{ij}\leq 1,
\label{cc-lin2}
\end{eqnarray}
where $s_{ij}$ is an auxiliary unconstrained real variable (used to handle the two original
constraints). Then, Eqs.~(\ref{cc-lin1},\ref{cc-lin2}),  which are linear in the optimization variables, as well as the respective linear inequalities originating from the generation CC are all added (accounted for) at any elementary step of the multi-step process, where an individual step consists in solving quadratic programming with linear constraints. As far as the nonlinear (but convex) constraint (\ref{cc-nonlinear}) is concerned,  we check if the constraint is valid at the current values of the optimization variables (known from the previous iteration). Valid constraints are ignored, while the violated ones are linearized around the current value, $\hat{\delta}$:
\begin{eqnarray}
&& C_{ij}(\delta)=\left[\sum_{k\in{\cal W}}
\sigma_k^2(\pi_{ik}-\pi_{jk}-\delta_i+\delta_j)^2\right]^{1/2}\leq s_{ij} \Rightarrow
\nonumber\\
&& C_{ij}(\hat{\delta})+\frac{\partial C_{ij}(\hat{\delta})}{\partial\delta_i}(\delta_i-\hat{\delta}_i)+\frac{\partial C_{ij}(\hat{\delta})}{\partial\delta_j}(\delta_j-\hat{\delta}_j)\leq s_{ij}.
\nonumber
\end{eqnarray}
A similar linearization of the active set of constraints is carried out with the generation-related CC violated constraints. At any new step all the nonlinear CC constraints are checked and the most violated is linearized and added to the active set.  The algorithm terminates after no violated (nonlinear) constraints are discovered. As discussed in the next Section,  only very few iterations (dozen or less) are required for typical CC-OPF case over an even very large network.

\section{VI. Empirical Results}

As discussed earlier, the inclusion of the sync constraints in our Sync CC-OPF increases the number of constraints that are checked during the execution of the cutting plane portion of the algorithm (the number of edge-related constraints is doubled).  However, in general this does not introduce any additional complexity issues beyond with what was discussed in \cite{12BCH}.

For completeness we will first recount the main contributions of our general approach (as first stated in Sections 2 and 3 of \cite{12BCH}).  We will then describe results that include the sync constraints.

\begin{itemize}
\item {\it Algorithm scalability.} \cite{12BCH} described results on large graphs including the various instances of the Polish grid model that consist of 2383-3120 edges, 327-388 generators and 2896-3693 lines and a BPA model with 2209 buses, 176 generators and 2866 lines.  Solving the CC-OPF required  5-30 seconds on a standard 4-core laptop. Even for these very large models the number of the cutting plane iterations was relatively small, i.e. 2-30, with only a handful of lines violating the chance constraints. Finally, in \cite{12BCH} the generation limits were sufficiently high such that that were non-binding.

\item {\it CC-OPF succeeds where standard OPF fails.} \cite{12BCH} shows that the standard OPF solution that ignores fluctuations in wind, may lead to choices of $p$ and $\alpha$ that overload lines too often. The CC-OPF optimization solutions  choose $p$ and $\alpha$  with significantly lower probability (risk) of line overloads, due to the utlization of the chance constrained parameters $\epsilon$.

\item  {\it Cost of Reliability.} From a certain point of view the CC-OPF produces solutions with lower cost than standard (i.e., not chance constrained) OPF.  In order for standard OPF to meet the same level of risk as CC-OPF, the amount of renewable energy must be reduced (when renewable energy has negligible generation cost), thereby driving generation costs higher.

\item  \cite{12BCH} also noted that there is no intuitively easy policy for ``fixing'' a standard OPF solution. In multiple experiments with different distributions of generation, line flows changed significantly in response to changes in the various governing parameters. In particular,  standard OPF and CC-OPF tend to produce very different flows patterns when the chance constraints are binding.

\item \cite{12BCH} also shows how the solution to CC-OPF helps when answering the following question:{\it What is the level of wind penetration that can be tolerated?} Given that increased wind penetration leads to increased risk of physical violations (e.g. of line overloads) a ``what-if'' experiment where wind farm output is progressively scaled up will find a feasible solution with the largest wind penetration possible. This critical feasible solution sets a threshold for reasonable investments in wind farms that are beneficial and do not require other upgrades of the grid.

\item  CC-OPF modeling also helps resolve other investment questions, for example the siting of wind farms. Different allocations of wind farm capacity over a grid typically result in very different solutions to CC-OPF and different risk exposure. In this case CC-OPF computation can be used as a diagnosis tool to identify nodes (or regions) in the grid where placement of wind-farms is desirable or prohibited.

\item Finally, the experiments of \cite{12BCH} have shown that (allowed) fluctuations in wind may introduce significantly variable operating conditions.  For example, flow reversals in response to minor changes in load, wind forecasts, and level of risk.

\end{itemize}

In the rest of this section we focus on analysis and illustrations that are specific to the synch constrained implementation of the CC-OPF, thereby increasing the contributions of CC OPF capabilities, advantages and accomplishments.

\begin{figure*}[t!] \centering
\includegraphics[width=1.5in]{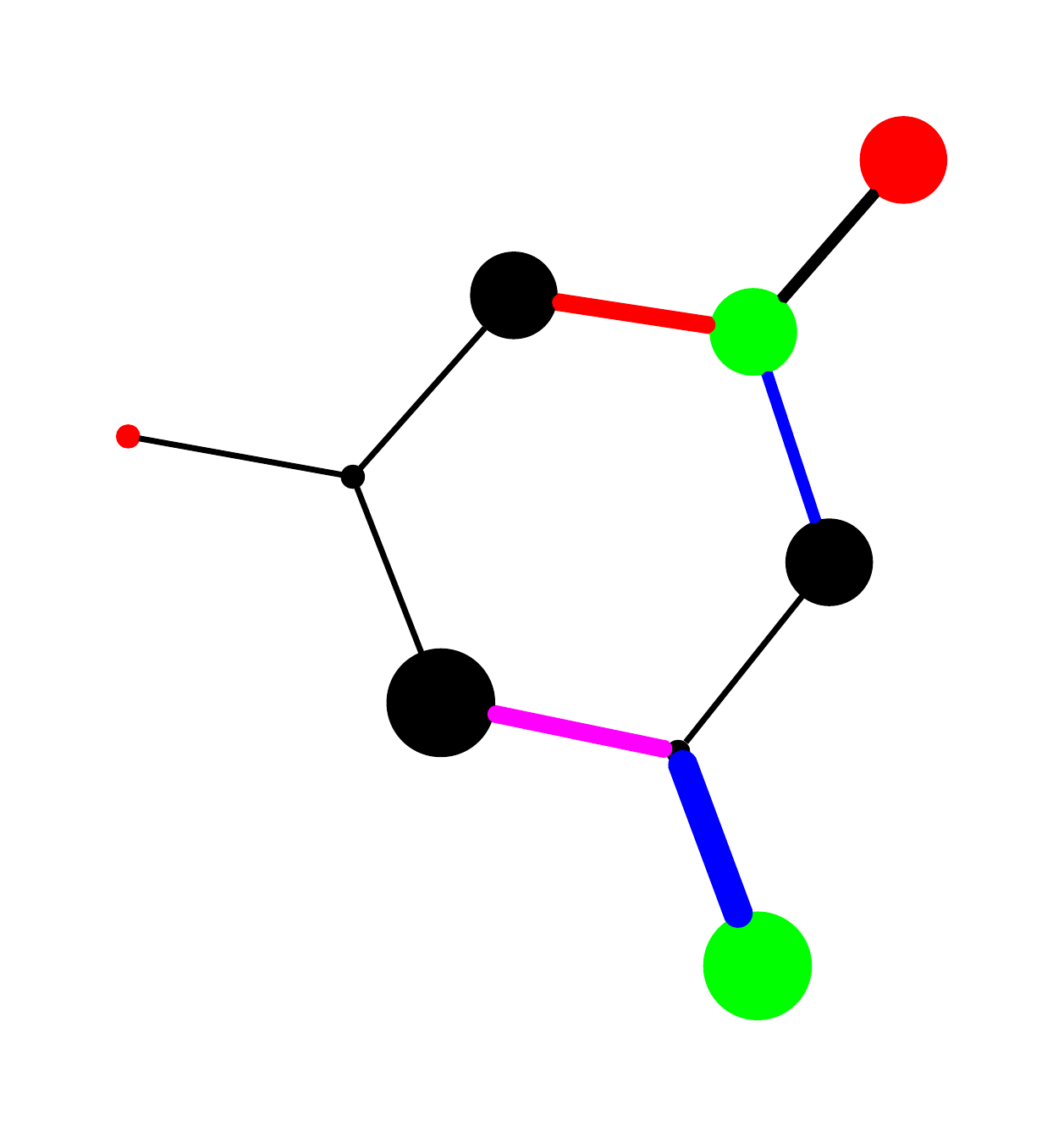}
\includegraphics[width=1.5in]{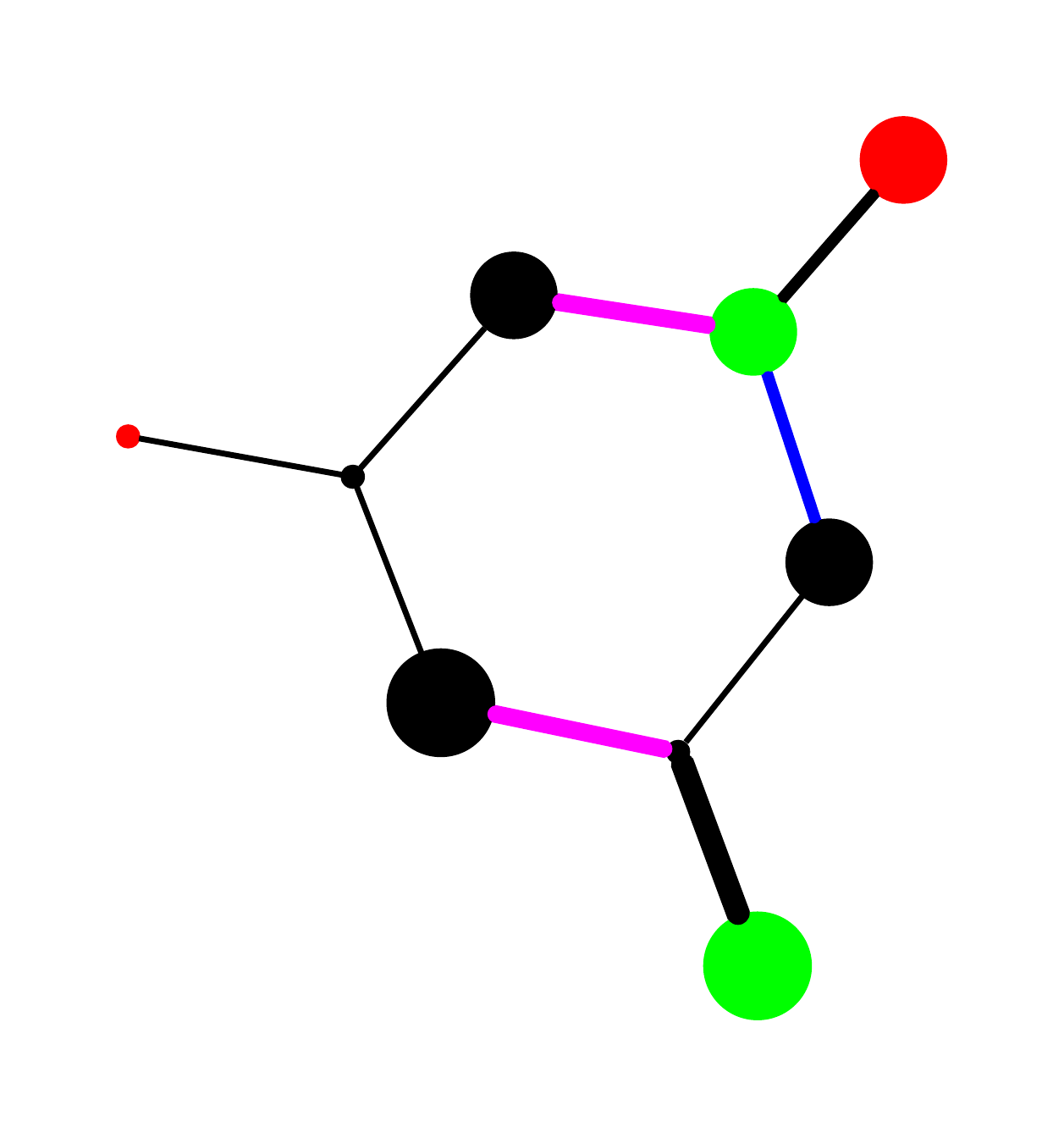}
\includegraphics[width=1.5in]{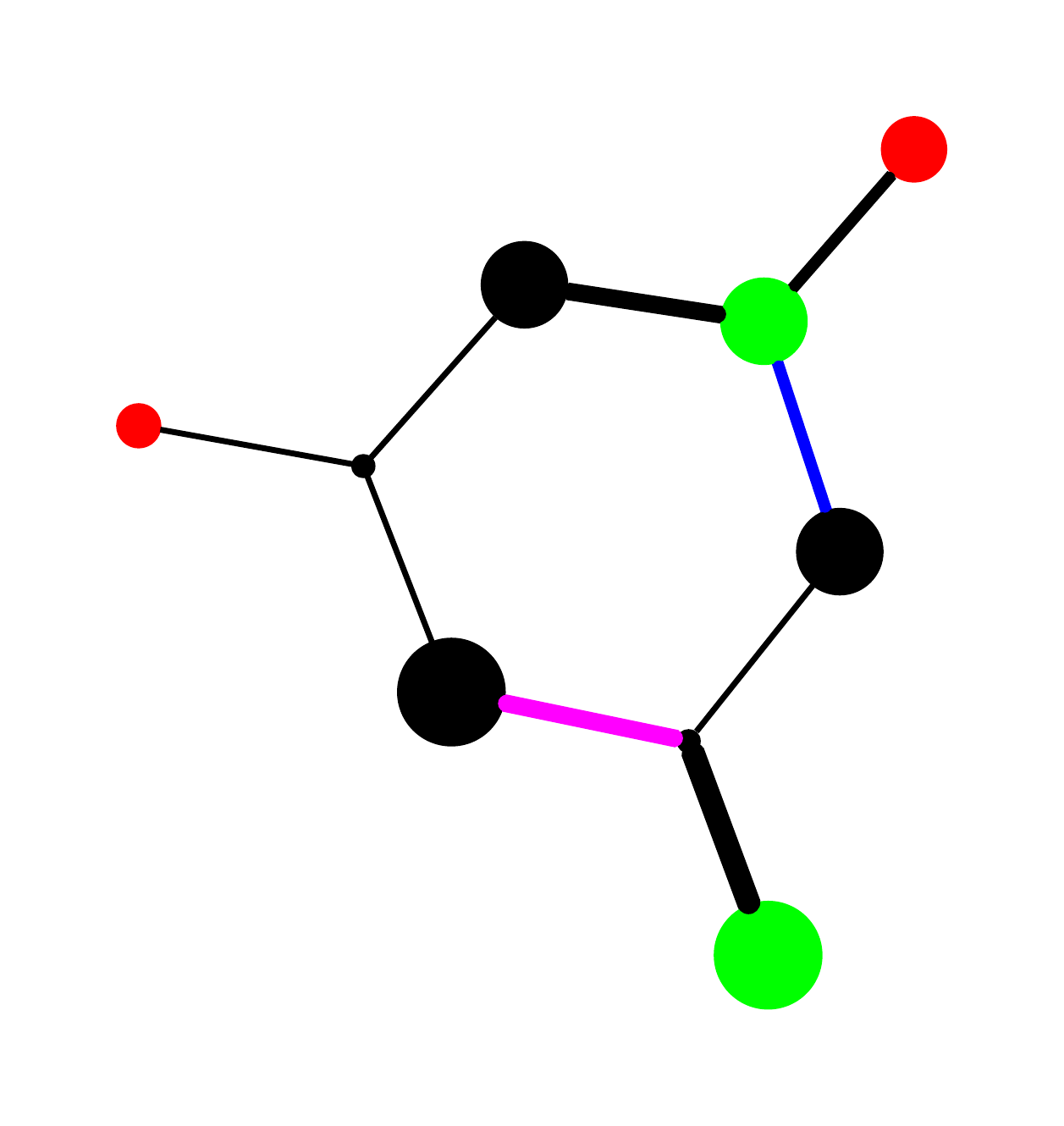}
\includegraphics[width=1.5in]{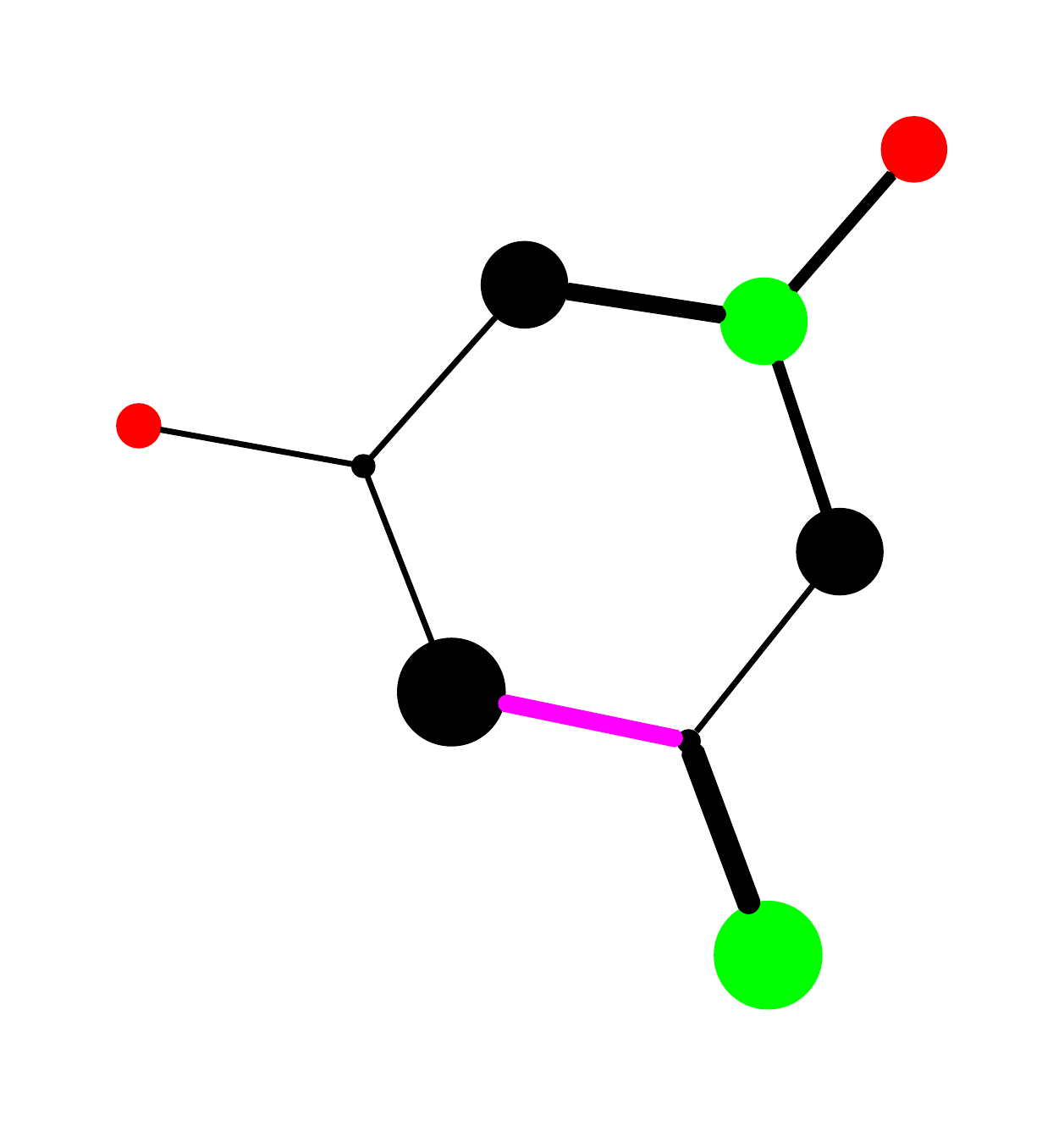}
\caption{
Outputs of the 1st, 8th, 11th and 13th iteration steps (in ascending order from left to right) of the cutting plane Sync CC-OPF algorithm for the 9 node IEEE model. The algorithm terminates after 13 iterations and these pictures show snapshots of the most significant qualitative changes. Loads, wind farms and regular generators are located at the nodes marked in black, green and red, respectively. The size of the nodes are scaled with consumption (of load), mean production (of wind farm) and optimal production (of regular generation). Red, magenta, blue and black show lines with both sync and thermal CC violations, only sync CC violations, only thermal CC violations, and no violations respectively.  The width of the lines are scaled according to the mean flow over the line.
\label{fig:9}} \end{figure*}

\subsection{Competition of sync and thermal risks guides iterations of the algorithm}

While at termination our algorithm produces a solution that satisfies all constraints, a useful empirical observation is that the algorithm ``discovers'' the set of lines that are most exposed to risk -- typically, these are the lines for which the conic constraints are violated during intermediate iterations of the algorithm, thus requiring the addition to cutting-planes as described above.    This phenomenon is easily explained by the (highly) nonlinear nature of the ``risk'' that is modeled by the chance constraints; this nonlinearity necessarily demands a comparatively larger number of inequalities so as to obtain an accurate approximation. Often, this set of critical lines is small and sometimes quite small.  There is, therefore, qualitative value in studying this set.

In this regard, an interesting pattern that emerged in our experiments was an alternation of violations in sync and thermal CC-constraints during the execution of the cutting plane algorithm.  This was observed in cases where $\bar{p}_{ij}/\beta_{ij}\leq 1$; even though the probabilistic (temporal) requirements on the loss of synchrony are much more stringent than the thermal line requirements,  i.e. $\varepsilon\ll \epsilon_{ij}$. Fig.~(\ref{fig:9}) illustrates this alternation of constraint violations. In this example CC constraints on four lines are violated after the first iteration. Of the four lines with CC violations, one line displays problems with both the thermal and sync conditions and one line has a violation in sync conditions only. The number of chance constraint violations is reduced as the number of iterations increases.  Naturally,  the optimal solution to generation dispatch changes during the course of the iterations.  In particular, we observe that as the algorithm iterates, flows over all lines with violated chance constraints are reduced. This is arguably driven by changes in generation dispatch over the grid and which is also accompanied by adjustment of the flow pattern over many other lines within the system.

\subsection{Pattern of Sync Warnings}


\begin{figure*}\centering
\includegraphics[width=2.5in]{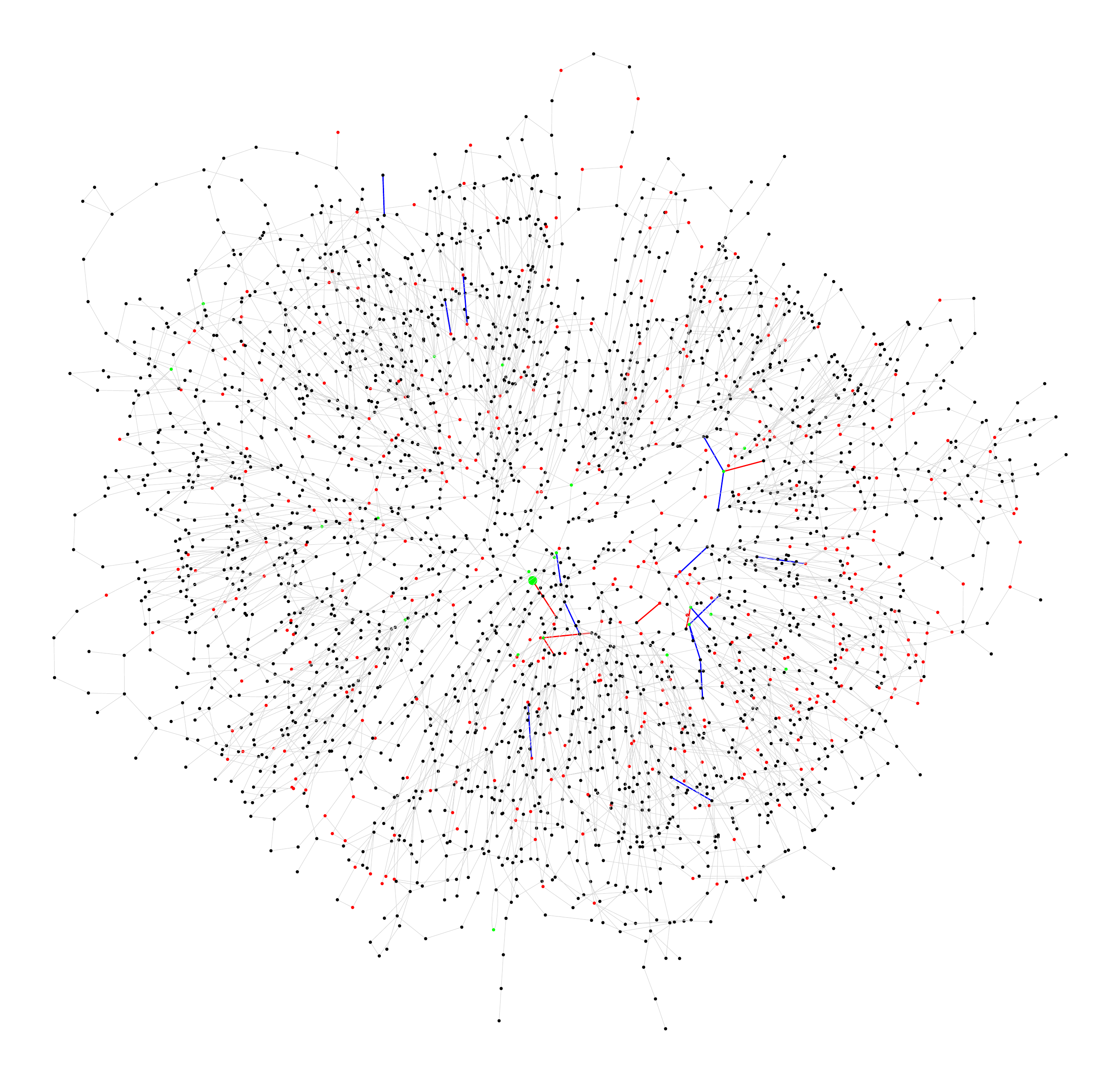}
\includegraphics[width=2.5in]{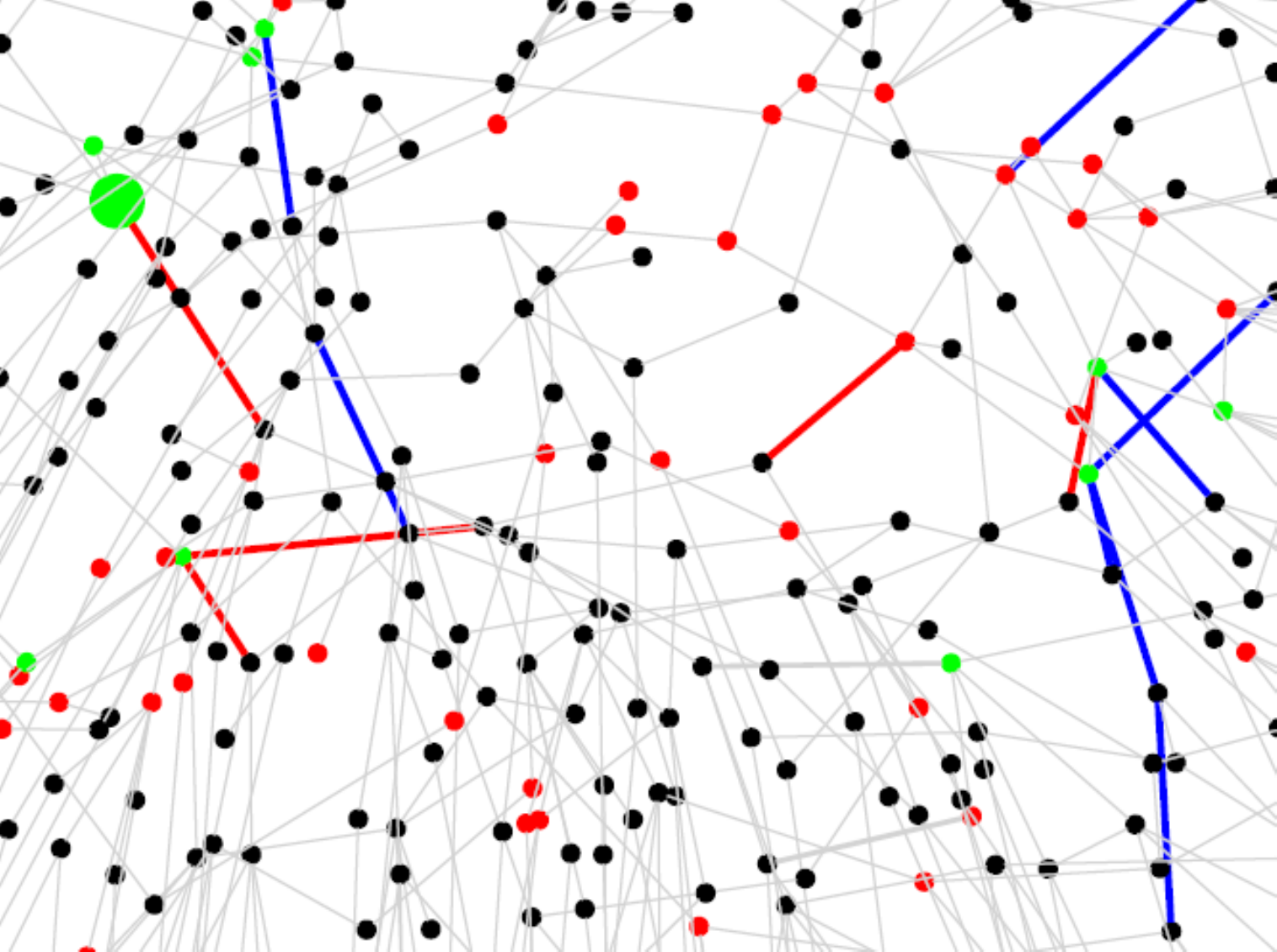}
\caption{A feasible solution snapshot of the Polish grid (rendered non-geographically) and a magnified spot from the snapshot are shown. Lines marked red and blue were sync-overloaded during the algorithm iterations.
Lines marked red showed the sync overload probability larger than $10^{-4}$ but smaller than $10^{-2}$ within the feasible solution,  where the latter was the pre-set value of $\varepsilon$. Scaling of node size and line width is set according to respective consumption/production and power flows within the final feasible solution.
\label{fig:Polish}} \end{figure*}


Fig.~(\ref{fig:Polish}) specifically focuses on sync chance constraints that are violated during the course of the algorithm; as explained above  the lines corresponding to these constraints are indicative of potentially vulnerable  patterns within the grid. In Fig.~(\ref{fig:Polish}) we color those lines that were overloaded at some iteration of the algorithm (which in this case terminated in 11 iterations). Lines marked red showed a sufficiently high probability of the sync overload even within the feasible solution (their sync probability of overload was larger than $10^{-4}$ but smaller than $\varepsilon=10^{-2}$).

\subsection{Sensitivity of the optimal solution to risk awareness and other parameters}

\begin{figure*}[t!] \centering
\includegraphics[width=2.5in]{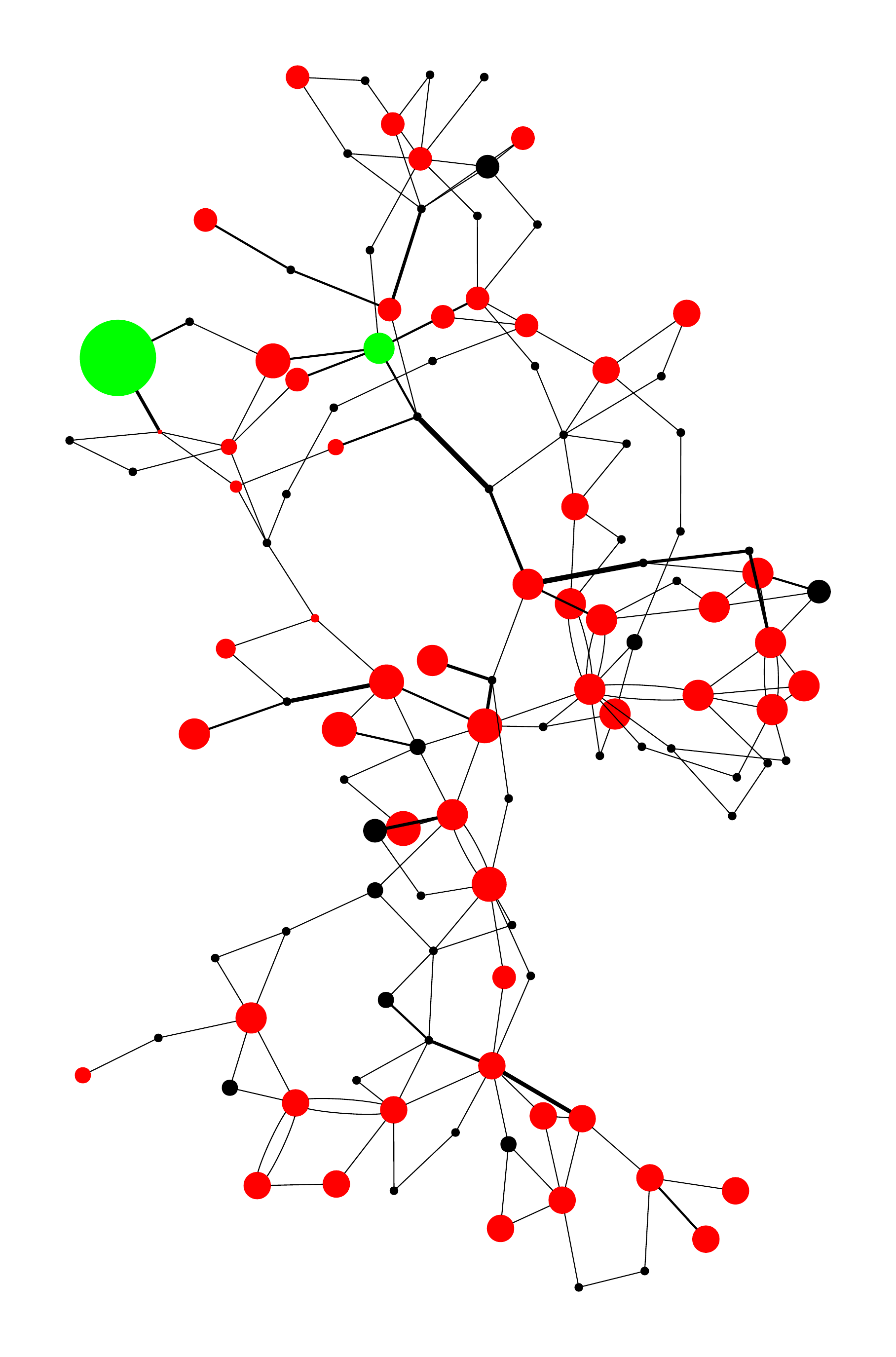}
\includegraphics[width=2.5in]{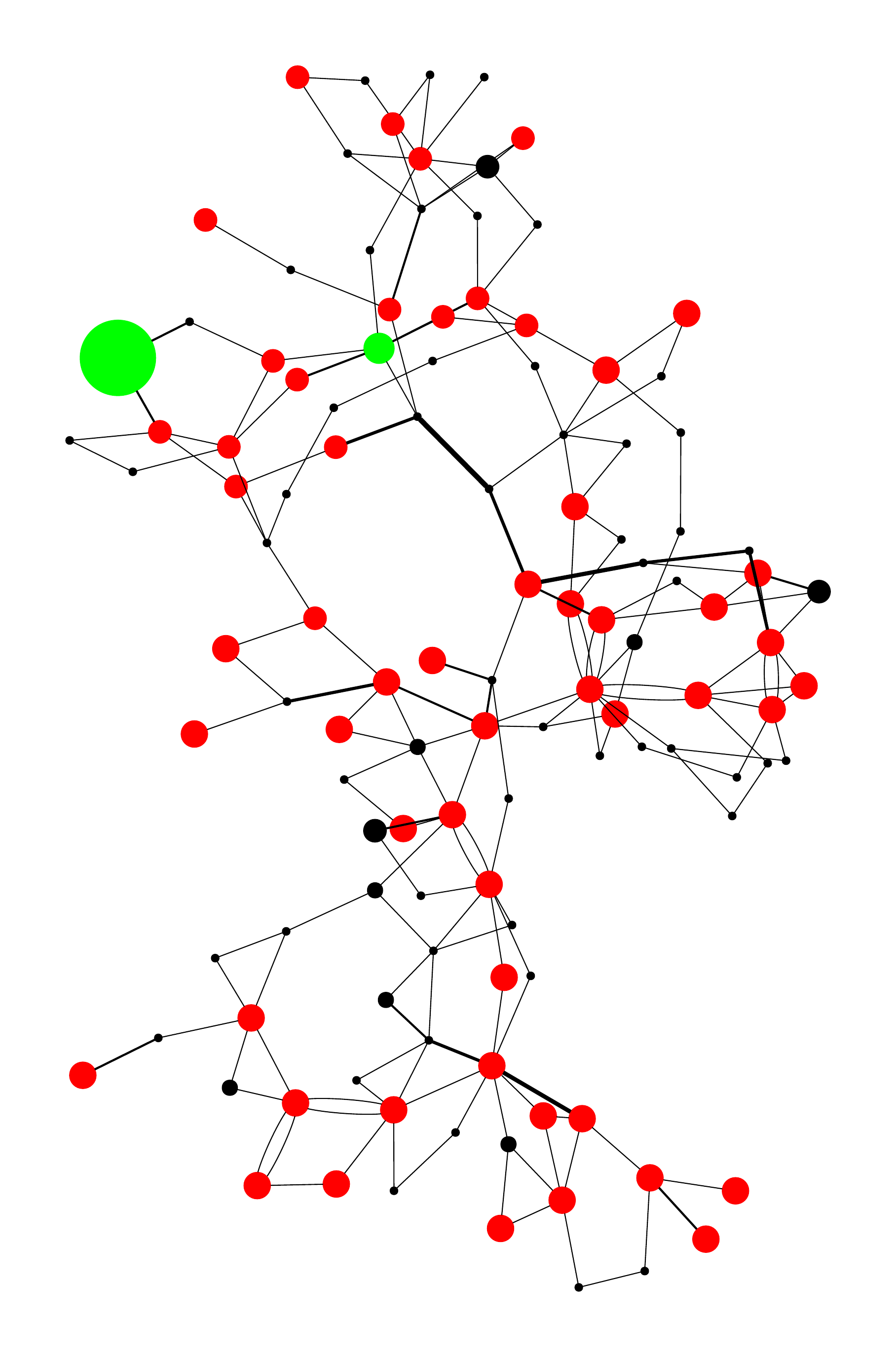}
\caption{
The two figures correspond to optimal solutions of the Sync CC-OPF on the 118 bus IEEE model. Green, black and red dots mark wind farms, loads and regular generators, respectively. The sizes are proportional to the mean production,  consumption and the average value of the optimal solution production. Line width expresses the mean value of power flows for the optimal solution. The two optimal solutions are derived for two slightly different values of the base voltage. The solution shown on the left required 21 iterations of the cutting plane algorithm with only sync conditions violated. The solution shown on the right required 9 iterations with violations  of both sync and thermal constraints. Optimal distribution of the mean generation is visibly less uniform for the sync constrained solution (shown on the left).
\label{fig:118}} \end{figure*}

The optimal solution to the Sync CC OPF problem may have significantly different structure depending on the governing parameters;  i.e. load, wind penetration and voltage level. A particularly interesting case, observed under two slightly different conditions, is illustrated in Fig.~(\ref{fig:118}) . In the case shown on the right the distribution of generation (red dots) is rather uniform,  while the same distribution is visibly much less uniform in the case shown on the left side of the figure. The left (non-uniform) case mainly has synchronization CC violations,  while the right (uniform) case contained a  mixture of the synchronization and thermal violations.  It is important to note that difference in solutions is genuine, as the two optimal values are different.

\section{VII. Conclusions and Path Forward}

This paper describes an approach for incorporating important nonlinear aspects of the AC power flow equations into OPF models. We first develop a formulation for AC-OPF that is convex under the assumption that (a) power lines are inductive (zero resistivity) and (b) voltage is maintained constant across the grid. This convex formulation builds on the earlier result of \cite{12BV}. Second, the paper developed a chance constrained AC-OPF that forces the probability of the grid losing synchrony,  the probability of thermal line overload and the probability of the generators deviating from their bounds to be sufficiently low. The AC OPF problem is reduced to a convex (conic) optimization problem that has the same complexity as the DC CC-OPF of \cite{12BCH}. Experiments show it to be efficiently solvable. Our formulation is approximate and it is based on implementing the static and linear synchronization conditions of \cite{12DCB}.

While these advances are considerable, there are a number of questions and challenges left for further explorations.  In particular,  we highlight the following open questions:
\begin{itemize}

\item Lemma 4.1 guarantees that the convex optimization (\ref{conserv2}) is a theoretically accurate approximation for AC-OPF when voltages are fixed and there are no losses. However, an efficient algorithm for solving this approximation is not yet developed. Moreover, inspired by \cite{12KCLB}, we believe that a distributed version of the algorithm is possible.

\item Adding chance constraints to the convex, nonlinear and exact AC-OPF makes finding a computationally efficient algorithm a  challenging task. One possible way of addressing this challenge, relies on approximating the chance constraints in a Large Deviation fashion and is briefly outlined in Appendix B.

\item As shown theoretically in \cite{12BCH} the robust version of the DC CC-OPF includes uncertainty in parameters of the wind forecast distribution. This approach generalizes to the sync CC-OPF introduced and discussed in this manuscript, but needs to be tested.

\item The convex AC-OPF and the sync CC-OPF do not include voltage variations and thus neglects dramatic and highly nonlinear phenomena such as voltage collapse \cite{68WC,75Ven,00Cut}. Including voltage-related effects is expected to be a much more difficult task than dealing with the phase angle nonlinearities and related loss of synchrony phenomena addressed in this paper.

\item Finally, there remain a large number of extensions to consider. These include modeling nonzero power line resistivity, including $N-1$-contingency compliance, modeling temporal (discrete or continuous) evolution of wind forecast and the incorporation of our approach into discrete problems such as unit commitment. These serve as a foundation for extending this  emerging approach to stochastic modeling and efficient algorithmic implementations of optimization and control problems in power transmission systems.

\end{itemize}

\section{Acknowledgments}

The work at LANL was carried out under the auspices of the National Nuclear Security Administration of the U.S. Department of Energy at Los Alamos National Laboratory under Contract No. DE-AC52-06NA25396. RB and MC also acknowledge partial support from LANL’s Grid Science project funded by the Advanced Grid Modeling program in U.S. DOE’s Office of Electricity Delivery and Energy Reliability and the DTRA Basic Research grant BRCALL08-Per3-D-2-0022. DB was partially supported by DOE award DE-SC0002676.

\appendix

{\bf \large A. Relation to Energy Function}

Eq.~(\ref{IndLoss}) can also be turned into an alternative formulation, where the only optimization is over $\theta$.  To derive such a formula, dual to Eq.~(\ref{IndLoss}), one starts from Eq.~(\ref{IndLoss}), incorporates the 
conditions (\ref{conserv}) through Lagrangian multipliers $\theta$ into the effective Lagrangian, performs variation over $\rho$, uses Eq.~(\ref{rho-theta}) to express $\rho$ via $\theta$, and finally one substitutes the result in the objective,  thus arriving at
\begin{eqnarray}
&&\!\!\!\!\max_{\theta} \!\left(\sum_{(i,j)\in{\cal E}} \beta_{ij}\left(\cos(\theta_i\!-\!\theta_j)\!-\!1\right)\!+\!2\sum_{i\in{\cal V}}\theta_i(p_i\!+\!d_i\!-w_i)\right)\nonumber\\
&&\!\!\!\!\!=\!2 \min_\theta \!\left(\sum_{(i,j)\in{\cal E}} \! \beta_{ij}\frac{1-\cos(\theta_i-\theta_j)}{2} -\sum_{i\in{\cal V}}\theta_i(p_i\!+\!d_i\!-\!w_i)\right).\nonumber
\end{eqnarray}
Note that the objective on the rhs of the last expression is nothing but the so-called energy function (modulo overall factor of $2$ and considered under fixed voltage) discussed extensively in the power systems literature, see e.g. \cite{89Pai} and references therein.

{\bf \large B. Large Deviation Derivation of the Nonlinear Chance-Constrained Condition}

Our goal for this Appendix is not to present the complete algorithm,  and even not to present a convex optimization problem,  but instead we aim to show how the CC constraints can be evaluated in the LD asymptotic,  thus resulting in some deterministic, but implicit and still non-convex expressions.  These expressions will require further manipulations/tricks to produce  an efficient computational procedure in the future.

Our task here is to replace the exact probabilistic chance-constrained condition of the thermal type, Eq.~(\ref{longchance-lin1}), e.g.
\begin{eqnarray}
\mbox{Prob}\left(\beta_{kl}\sin(\theta_k-\theta_l)\geq \varrho\right) < \epsilon,
\label{CC1}
\end{eqnarray}
by its Large-Deviation (LD), or saddle-point (instanton) deterministic approximation.
Then, using Eq.~(\ref{Gauss}) for the PDF of wind one approximates the rhs of Eq.~(\ref{CC1}) in the LD (saddle-point) manner
\end{multicols}
\begin{eqnarray}
&& \mbox{Prob}\left(\beta_{kl}\sin(\theta_k-\theta_l)\geq \varrho\right)=
\int\limits_{\omega\in \Omega(\bar{p},\alpha)}
  d\omega \rho(\omega)= W_0 \int\limits_{\omega\in \Omega_0(\bar{p},\alpha)}
  d\theta  \rho(\omega)=W_1 \exp\left(-E(\bar{p};\alpha)\right),\label{LD}\\ &&
\Omega(\bar{p},\alpha)=\left(\omega\in \mathbb{R}^{|{\cal W}|}\left|\ \exists \theta\in \mathbb{R}^{|{\cal V}|} \mbox{ s.t. }
\left\{\begin{array}{c}
\theta_n=0\\ \forall i\in{\cal V}:\
\sum_{j:(i,j)\in{\cal E}}\beta_{ij}\sin(\theta_i-\theta_j)=(\bar{p}-d+\mu+\omega-(e^T\omega)\alpha)_i\\
\beta_{kl}\sin(\theta_k-\theta_l)\geq\varrho
\end{array}\right.\right.\right),
\label{Omega}\\
&& \Omega_0(\bar{p},\alpha)=\left(\omega\in \mathbb{R}^{|{\cal W}|}\left|\ \exists \theta\in \mathbb{R}^{|{\cal V}|} \mbox{ s.t. }
\left\{\begin{array}{c}
\theta_n=0\\ \forall i\in{\cal V}:\
\sum_{j:(i,j)\in{\cal E}}\beta_{ij}\sin(\theta_i-\theta_j)=(\bar{p}-d+\mu+\omega-(e^T\omega)\alpha)_i\\
\beta_{kl}\sin(\theta_k-\theta_l)=\varrho
\end{array}\right.\right.\right),
\label{Omega0}\\
&& E(\bar{p},\alpha)=\left.\min_{\theta,\vartheta,\omega} \sum_{i\in{\cal W}}\frac{ \omega_i^2}{2 \sigma_i^2}\right|_{\begin{array}{c}\beta_{kl}\vartheta_{kl}=\varrho,\quad \theta_n=0,\\
\forall i\in{\cal V}:\quad \sum_{j:(i,j)\in{\cal E}}\beta_{ij}\vartheta_{ij}=(\bar{p}-d+\mu+\omega-(e^T\omega)\alpha)_i,\\
\forall (i,j)\in{\cal E}:\quad\vartheta_{ij}=-\vartheta_{ij}=\sin(\theta_i-\theta_j)
\end{array}},
\label{E}
\end{eqnarray}
\begin{multicols}{2}
where $W_0$ and $W_1$ are volume factors which dependence on $\sigma$ is expected to be algebraic in the $|\sigma|\to 0$ limit. Then, the condition (\ref{CC1}) turns into
\begin{equation}
E(\bar{p};\alpha)\geq \log(1/\varepsilon)+\log W_1.
\label{E_cond}
\end{equation}
The ``effective energy", $E(\bar{p};\alpha)$, dependence on $|\sigma|$ is expected to be at least $\sim 1/|\sigma|$ at $|\sigma|\to 0$ (or stronger), therefore one can safely ignore the $\log W_1=O(\log\sigma)$ term on the rhs of Eq.~(\ref{E_cond}) in the limit.

As shown in the Subsection below the linear (DC) approximation version of Eq.~(\ref{E_cond}) results in a constraint which is convex in $\bar{p}$ and $\alpha$. An important question becomes: if Eq.~(\ref{E_cond}) may also allow a convex and/or computationally convenient re-formulation?

\subsection{The case of DC approximation}

Eq.~(\ref{LD}) applies to the nonlinear case of interest (fixed voltage, no resistive losses),  however under a simple modification,  consisting in replacement of $\vartheta_{ij}\to (\theta_i-\theta_j)$, it reduces to
\end{multicols}
\begin{equation}
E_{\mbox{DC}}(\bar{p};\alpha)=\left.\min_{\theta,\omega} \sum_{i\in{\cal W}}\frac{ \omega_i^2}{2 \sigma_i^2}\right|_{\begin{array}{c}\beta_{kl}(\theta_k-\theta_l)=\varrho,\quad \theta_n=0,\\
\forall i\in{\cal V}:\quad \sum_{j:(i,j)\in{\cal E}}\beta_{ij}(\theta_i-\theta_j)=(\bar{p}-d+\mu+\omega-(e^T\omega)\alpha)_i\end{array}},
\label{E-DC}
\end{equation}
\begin{multicols}{2}
describing the CC-constraint within the DC-approximation.  Then the DC analog of
Eq.~(\ref{E_cond}) becomes
\begin{equation}
E_{\mbox{DC}}(\bar{p};\alpha)\geq \log(1/\varepsilon)\pm O(\log \sigma).
\label{E_cond_DC}
\end{equation}

Resolving the conditions in Eq.~(\ref{E-DC}) explicitly and following the notations of \cite{12BCH} one derives
\end{multicols}
\begin{eqnarray}
&& E_{\mbox{DC}}(\bar{p};\alpha)=\left.\min_{\omega} \sum_{i\in{\cal W}}\frac{ \omega_i^2}{2 \sigma_i^2}\right|_{
\left(\breve{B}(\bar{p}-d+\mu+\omega-(e^T\omega)\alpha)\right)_k-
\left(\breve{B}(\bar{p}-d+\mu+\omega-(e^T\omega)\alpha)\right)_l=\varrho}
\label{DC2}\\
&&=\min_{\omega} \max_\phi \left(\sum_{i\in{\cal W}}\frac{ \omega_i^2}{2 \sigma_i^2}-\phi\left(\left(\breve{B}(\bar{p}-d+\mu+\omega-(e^T\omega)\alpha)\right)_k
-\left(\breve{B}(\bar{p}-d+\mu+\omega-(e^T\omega)\alpha)\right)_l-\varrho\right)\right)
\label{DC3}\\
&& =-\max_\phi \left(\phi^2\sum\limits_{i\in{\cal W}}\frac{\sigma_i^2}{2}
\left(\sum\limits_{j\in{\cal V}}(\breve{B}_{kj}-\breve{B}_{lj})\alpha_j-\breve{B}_{ki}+\breve{B}_{li}\right)^2\!+\!
\phi\left(\left(\breve{B}(\bar{p}-d+\mu)\right)_k
\!-\!\left(\breve{B}(\bar{p}-d+\mu)\right)_l\!-\!\varrho\right)\right)
\label{DC4}\\
&& =\frac{\left(\left(\breve{B}(\bar{p}-d+\mu)\right)_k
\!-\!\left(\breve{B}(\bar{p}-d+\mu)\right)_l\!-\!\varrho\right)^2}{2
\sum\limits_{i\in{\cal W}}\sigma_i^2
\left(\sum\limits_{j\in{\cal V}}(\breve{B}_{kj}-\breve{B}_{lj})\alpha_j-\breve{B}_{ki}+\breve{B}_{li}\right)^2},
\label{DC5}
\end{eqnarray}
\begin{multicols}{2}
where the auxiliary variational variable $\phi$ is a scalar. Eq.~(\ref{DC5}) combined with Eq.~(\ref{E_cond_DC}) results in the $\sigma_i\to 0$ (and $\beta_{ij}\to 1$) version of the CC-constraint (\ref{longchance-lin1}), with
$\eta(x)\to (2\log(1/x))^{1/2}$.

\bibliographystyle{siam}
\bibliography{Bib/abstract}

\end{multicols}

\end{document}